\documentclass[12pt]{article}

\usepackage{amssymb}

\usepackage{amsthm}

\title{Borel transforms of functions from a parameterized family of Hilbert spaces of entire functions}

\author{K.P. Isaev, R.S. Yulmukhametov}

\begin{document}

\maketitle

\section{Introduction}
\label{Intr}

Let $D$ be a bounded convex domain on the complex plane. The Bergman space $B_2(D)$ and the Smirnov space $E_2(D)$ are quite well studied due to their importance in complex analysis problems. The space $B_2(D)=H(D)\bigcap L_2(D)$, where $H(D)$ is the space of analytic functions on $D$ and $L_2(D)$ is the space of functions with square-summable modulus. The space $E_2(D)$ consists of functions $f\in H(D)$ having boundary values of the class $L_2(\partial D)$. In [1], [2], in particular, it was found that the Laplace transform of linear continuous functionals $\mathcal L:S \longrightarrow \widehat S(\lambda )=S_z(e^{\lambda z})$ establishes an isomorphism between the space dual to $E_2(D)$ and the space of entire functions
$$
\widehat E_2(D) =\left \{ F\in H(\mathbb C): \int  _0^{2\pi }\int _0^\infty \frac {|F(re^{i\varphi })|^2drd\Delta (\varphi  )}{K_1(re^{i\varphi })}<\infty \right \},
$$
where $h(\varphi )$ the support function of the domain $D$, that is,
\begin{eqnarray*}
& h(\varphi)=\max _{z\in \overline D} \textrm {Re } ze^{i\varphi },\quad \varphi \in [0;2\pi ],\\
& \Delta (\varphi )=h(\varphi )+\int _0^\varphi h(\theta )d\theta ,\quad \varphi \in [0;2\pi ],\\
& K_1(\lambda )=\|e^{\lambda z}\|^2_{L_2(\partial D)}=\int _D|e^{\lambda z}|^2ds(z),\quad \lambda \in \mathbb C,
\end{eqnarray*}
and $ds(z)$ is an element of length along the boundary of $D$.
In [3] it is shown that the Laplace transform establishes a continuous linear isomorphism between the space dual to $B_2(D)$ and the space of entire functions
$$
\widehat B_2(D) =\left \{ F\in H(\mathbb C): \int  _0^{2\pi }\int _0^\infty \frac {|F(re^{i\varphi })|^2drd\Delta (\varphi  )}{K(re^{i\varphi })}<\infty \right \},
$$
where
$$
K(\lambda )=\|e^{\lambda z}\|^2_{L_2(D)}=\int _D|e^{\lambda z}|^2dm(z),\quad \lambda \in \mathbb C.
$$
Considering that $K_1(\lambda )\cong |\lambda |K(\lambda )$, $|\lambda |\longrightarrow \infty $, for $\beta \in \mathbb R$ it is natural to introduce the spaces
$$
P_\beta (D)=\left \{ F\in H(\mathbb C):\ \|F\|^2:=\int _0^{2\pi }\int _0^\infty \frac {|F(re^{i\varphi })|^2drd\Delta (\varphi  )}{K(re^{i\varphi })r^{2\beta }}<\infty \right \}.
$$
The spaces $P_\beta (D)$ form a continuous scale of Hilbert spaces, and, as stated above, $\widehat P_0(D)$ is isomorphic to $B_2(D)$ and $\widehat P_{1/2} (D)$ is isomorphic $E_2(D)$. Thus, Bergman and Smirnov spaces are embedded in the scale of Hilbert spaces, which, according to the authors, should make it possible to apply the theory of Hilbert scales to study of problems in these spaces.

The function associated with an entire function $F$ of exponential type $\sigma $ in the sense of Borel is
$$
\gamma (z)=\sum _{k=0}^\infty \frac {F^{(k)}(0)}{z^{k+1}}, \quad |z|>\sigma .
$$
Let $H_0(\mathbb C\setminus \overline D)$ be the space of functions that are analytic in $\mathbb C\setminus \overline D$ and vanish at infinity. Let $G^{(\alpha )}(D)$ denote the space
$$
\left \{h(\zeta )\in H_0(\mathbb C\setminus \overline D):\ \|h\|^2:=\int _{\mathbb C\setminus \overline D}|h''(\zeta )|^2 {\textrm { dist } }^{2\alpha }(D,\zeta )\, dm(\zeta )<\infty \right \}
$$
for $\alpha >0$. In this paper we intend to prove the following theorem.

\noindent {\bf Main theorem.} {\it Let $F$ be an entire function with indicator diagram $D$, $\gamma $ be the function associated with $F$ in the sense of Borel, and $\beta \in (-\frac 12;\frac 32)$.
Then for some constants $c(\beta ,D), C(\beta ,D)>0$, depending on the domain $D$ and the parameter $\beta $, the following relation holds:
$$
c(\beta ,D)\|\gamma \|^2_{G^{\beta +1} }\le \|F\|_{P_\beta } =\int _0^\infty \int _0^{2\pi }\frac {|F(re^{i\varphi })|^2}{K(re^{i\varphi })r^{2\beta }}d\Delta (\varphi )dr\le C(\beta ,D)\|\gamma \|^2_{G^{\beta +1} }.
$$}

The interest in these spaces is caused by the fact that the authors have a reasonable assumption that the spaces $P_\beta(D)$ admit unconditional bases of reproducing kernels (see [5] ---[9]). Accordingly, a scale of Hilbert spaces of functions that are analytic in the convex domain $D$ and admit unconditional bases of exponentials will be found.

In Sections 2 and 3 we will prove the preparatory theorems. In Section 2 we will estimate the integral
$$
\int _0^\infty \frac {|F(re^{i\varphi })|^2dr}{K(re^{i\varphi })r^{2\beta }}
$$
for $\beta >-\frac 12$. In Section 3, the localization of the norm in the spaces $G^\alpha$ will be justified.

\section{Estimation of the integral with respect to the radius}
\label{Prel1}

For a fixed $\varphi \in [0;2\pi ]$, using the map $z\longrightarrow w=ze^{i\varphi }-h(\varphi )$, we transform $D$ into a domain $D_\varphi $ which lies in the left half-plane and is tangent to the axis of ordinates. For $t<0$ we denote by $s(t,\varphi )$ the area of the intersection of the domain $D_\varphi $ and the strip $\{z=x+iy:\ t<x<0 \}$.

\noindent {\bf Theorem 1}  {\it Let $F$ be an entire function satisfying the condition: for some $\beta \in (-\frac 1 2 ;+\infty )$ and $\varphi \in [0;2\pi ]$
$$
I_\varphi =\int _0^\infty \frac {|F(re^{i\varphi })|^2}{K(re^{i\varphi })r^{2\beta }}dr <\infty,
$$
and $\gamma $ is the Borel transform of $F$. Then the inequalities
$$
a(\beta )I_\varphi \le \int _{-\infty }^\infty \int _{-\infty }^0\frac {|\gamma ''(e^{-i\varphi }(h(\varphi )-(x+iy))|^2 |x|^{2\beta +3} }{s(x,\varphi )}dxdy\le A(\beta )I_\varphi
$$
hold for some constants $a(\beta ), A(\beta )$ depending only on the parameter $\beta $ (see the remark at the end of the section).
}
\begin{proof}
The proof is essentially based on the reasoning of $\S1$ in [3]. Therefore, we tried to maintain the appropriate designations. In the half-plane
$P_\varphi =\{ \zeta :\ \textrm {Re } \zeta e^{i\varphi }>h(\varphi )\}$ the function $\gamma $ is restored by the formula
$$
\gamma (\zeta )=\int _0^\infty F(re^{i\varphi })e^{-\zeta re^{i\varphi }}e^{i\varphi }dr.
$$
We write $\zeta =(h(\varphi )-\xi)e^{-i\varphi }$ for points in this half-plane, where $\xi $ varies over the left half-plane $\textrm {Re } \xi <0$. Then
\begin{equation}\label{1}
\gamma ''(e^{-i\varphi }(h(\varphi )-\xi))=\int _0^\infty \left (F(re^{i\varphi })r^2e^{3i\varphi }e^{-h(\varphi )r }\right )e^{\xi r}dr.
\end{equation}
Let's use the theorem from [4].

\noindent {\bf Theorem A}  {\it Let $u(t)$ be a convex function on an interval $I$ and let $\widetilde u$ be the Young conjugate of $u$,
$$
\widetilde u(x)=\sup _{t\in I}(xt-v(t)).
$$
We set
$$
K_0(x)=\int _Ie^{2xt-2v(t)}dt, \quad x\in J=\{x\in \mathbb R:\ \widetilde v(x)<\infty \}.
$$
Then, for an arbitrary function $g$ on $I$ such that the integral
$$
\|g\|^2=\int _I|g(t)|^2e^{-2v(t)}dt,
$$
is finite, the function $\widehat g(z)=\int _I\overline g(t)e^{zt-2v(t)}dt$ satisfies the inequalities
$$
a\|g\|^2\le \int _{-\infty }^{\infty }\int _J\frac {\widehat g(x+iy)|^2}{K_0(x)}d\widetilde v'(x)dy\le A\|g\|^2,
$$
where $a$ and $A$ are absolute positive constants, that is, they do not depend on the function $g$ or the weight $v$.
}

Let
$$
\eta (r)=\frac {e^{2h(\varphi )r}}{K(re^{i\varphi })},\quad u(r)=\frac 1 2 \ln \frac {\eta  (r)}{r^4},\quad v(r)=u(r)-\beta \ln r.
$$
Using these functions, we write formula (1) in the form
$$
\gamma ''(e^{-i\varphi }(h(\varphi )-\xi))=\int _0^\infty \frac {F(re^{i\varphi})e^{3i\varphi}e^{h(\varphi )r}}{r^{2+2\beta}K(re^{i\varphi})}e^{\xi r-2v(r)}dr.
$$
Next we intend to apply Theorem A, assuming $I=(0;+\infty)$,
$$
g(r)=\frac {\overline F(re^{i\varphi})e^{-3i\varphi}e^{h(\varphi )r}}{r^{2+2\beta}K(re^{i\varphi})},
$$
and taking $v(t)$ as the weight function. To do this, make sure that the function $v(t)$ is convex.

\noindent {\bf Lemma 1}

{\it \noindent 1) The function $v$ is convex on $I=(0;+\infty  )$, and
$$
\left (\frac 1 2+\beta \right )\frac 1{r^2}\le v''(r)\le (2+\beta )\frac 1{r^2},\quad r>0;
$$
\noindent 2) If $\widetilde v(x)=\sup _{r> 0} (xr-v(r))$, $x\in J=(-\infty ;0)$, is the Young conjugate of $v$, then
$$
\frac {(1+2\beta )^2}{4(2+\beta )}\frac 1{x^2}\le \widetilde v''(x)\le \frac {2(2+\beta )^2}{1+2\beta}\frac 1{x^2},\quad x<0.
$$
}
\begin{proof}
In the first item of Lemma 1 in [3], the following inequalities are proved
$$
\frac 1{2r^2}\le u''(r)\le \frac 2{r^2},\quad r>0.
$$
The estimates of item 1 follow from these inequalities. Let's prove the second item. The proof of Lemma 1 in [3] shows (see (5)) that
$$
-\frac 2{r}\le u'(r)\le -\frac 1{2r},\quad r>0.
$$
Therefore
\begin{equation}\label{2}
-( 2+\beta )\frac 1{r}\le v'(r)\le -\left (\frac 1 2+\beta \right )\frac 1 r.
\end{equation}
The following relations are also proved there
$$
u(r)\longrightarrow -\infty ,\quad \left |\frac {u(r)}{\ln r}\right |=O(1), \quad r\longrightarrow +\infty .
\quad \lim _{r\longrightarrow 0+}u(r)=+\infty ,
$$
Hence, the function $\widetilde v(x)$ is defined on $(-\infty ;0)$ and
$$
\lim _{r\longrightarrow 0+}(xr-v(r))=\lim _{r\longrightarrow +\infty }(xr-v(r))=-\infty .
$$
Thus, the supremum in the definition of $\widetilde v$ is attained at the single point $r=r(x)>0$ for which $v'(r)=x$. From estimates (2) we obtain
\begin{equation}\label{3}
-\left (\frac 1 2+\beta \right )\frac 1 x\le r(x)\le -( 2+\beta )\frac 1{x}.
\end{equation}
By the definition of $\widetilde v(x)$, we have the identity
$$
\widetilde v(x)\equiv xr(x)-v(r(x)),\quad x<0,
$$
or
$$
\widetilde v(v'(r))\equiv v'(r)r-u(r),\quad r>0.
$$
We differentiate the above identity twice:
$$
\widetilde v''(v'(r))v''(r)\equiv 1, \quad r>0.
$$
Hence, in view of formula (3) and assertion 1) of this Lemma, it follows that
$$
\frac {(1+2\beta )^2}{4(2+\beta )}\frac 1{x^2}\le \widetilde v''(x)=\frac 1{v''(r(x))}\le \frac {2(2+\beta )^2}{1+2\beta}\frac 1{x^2},\quad x<0.
$$
Lemma 1 is proved.
\end{proof}

We apply Theorem A to the function
$$
g(r)=\frac {\overline F(re^{i\varphi})e^{-3i\varphi}e^{h(\varphi )r}}{r^{2+2\beta}K(re^{i\varphi})},
$$
which shows that the expression
$$
\|g\|^2=\int _0^\infty |g(r)|^2e^{-2v(r)}dr=\int _0^\infty \frac {|F(re^{i\varphi })|^2dr}{K(re^{i\varphi })r^{2\beta }}=I_\varphi
$$
is equivalent to the integral
$$
\int _{-\infty }^\infty \int _{-\infty }^0\frac {|\gamma ''(e^{-i\varphi }(h(\varphi )-(x+iy)))|^2}{K_0(x,\varphi )}d\widetilde v'(x)dy.
$$
Estimating $d\widetilde v'(x)$ according to item 2 of Lemma 1, we obtain
$$
\frac {(1+2\beta )^2}{4(2+\beta )}aI_\varphi \le \int _{-\infty }^\infty \int _{-\infty }^0\frac {|\gamma ''(e^{-i\varphi }(h(\varphi )-(x+iy)))|^2}{K_0(x,\varphi )x^2}d\widetilde v'(x)dy\le
$$
\begin{equation}\label{4}
\le \frac {2(2+\beta )^2}{1+2\beta}AI_\varphi .
\end{equation}

To complete the proof of Theorem 1, we need one more lemma.

\noindent {\bf Lemma 2} {\it For arbitrary $\varphi \in [0;2\pi ]$ and $\beta \in (-\frac 1 2;\infty )$, the inequalities
$$
2^{-(2\beta +5)}a_0(\beta )\frac {s(t,\varphi )}{|t|^{2\beta +5}}\le K_0(t,\varphi )\le a_0(\beta )\left (1+\frac {a_+(\beta )}{a_-(\beta )}\right )\frac {s(t,\varphi )}{|t|^{2\beta +5}},
$$
hold, where
$$
a_0(\beta )=\int _0^\infty t^{2\beta +4}e^{-2t}dt,\quad a_-(\beta )=\int _0^1\frac {tdt}{(1+t)^{2\beta +5}},
$$
$$
a_+(\beta )=\int _1^\infty \frac {tdt}{(1+t)^{2\beta +5}}.
$$
}
\begin{proof}
Let
$$
D_\varphi =\{z=x+iy:\ f_1(x)<y<f_2(x),\quad R_\varphi <x<0\}.
$$
Then $f(x)=f_1(x)-f_2(x)$ is a non-negative concave function on the interval $[-R_\varphi ;0]$ and
$$
K_0(t,\varphi )=\int _0^\infty e^{2tr-2v(r)}dr=\int _0^\infty e^{2rt}\frac {r^{2\beta +4}dr}{\eta (r)}=
$$
$$
=\int _0^\infty e^{2xr}r^{2\beta +4}K(re^{i\varphi})e^{-2rh(\varphi )}dr=\int _0^\infty e^{2rt}r^{2\beta +4}\left (\int _{D_\varphi }e^{2rx}dxdy\right )dr=
$$
$$
=\int _{D_\varphi }\left (\int _0^\infty e^{2r(x+t)}r^{2\beta +4}dr\right )dxdy\quad \textrm{for}\quad t<0.
$$
Therefore
$$
K_0(t,\varphi )=a_0(\beta )\int _{D_\varphi }\frac {dxdy}{|x+t|^{2\beta +5}},
$$
where
$$
a_0(\beta )=\int _0^\infty e^{-2\tau }\tau ^{2\beta +4}d\tau.
$$

1.  Let $t\le -D_\varphi $, then $|t|\le |t+x|\le 2|t|$ for values of $x$ in the interval of integration,
and therefore
$$
\frac {a_0(\beta )2^{-(2\beta +5)}}{t^{2\beta +5}}|D_\varphi |\le K_0(t,\varphi )\le \frac {a_0(\beta )}{t^{2\beta +5}}|D_\varphi |,
$$
where $|D_\varphi |$ is the area of $D_\varphi $, which means that the assertion of the lemma is true since $s(t,\varphi )=|D_\varphi |$.

2. Let $0\ge t> -D_\varphi $ and $p=f(t)$. The concavity of $f$ implies the inequalities
$$
f(x)\le \frac p t x,\quad -R_\varphi \le x\le t,\quad f(x)\ge \frac p t x,\quad t\le x\le 0.
$$
Therefore
$$
\int _{-R_\varphi }^t\frac {f(x)dx}{|t+x|^{2\beta +5}}\le \frac p{|t|}\int _{|t|}^\infty \frac {rdr}{(r+|t|)^{2\beta +5}}=\frac { a_+(\beta  )p}{|t|^{2\beta +4}},
$$
$$
\int _{t}^0\frac {f(x)dx}{|t+x|^{2\beta +5}}\ge \frac p{|t|}\int _{0}^{|t|} \frac {rdr}{(r+|t|)^{2\beta +5}}=\frac { a_-(\beta  )p}{|t|^{2\beta +4}},
$$
where
$$
a_+(\beta )=\int _1^{\infty }\frac {\tau d\tau }{(1+\tau )^{2\beta +5}},\quad a_-(\beta )=\int _0^{1}\frac {\tau d\tau }{(1+\tau )^{2\beta +5}},
$$
which means that
$$
K_0(t,\varphi )=a_0(\beta )\int _{-R_\varphi }^0\frac {f(x)dx}{|t+x|^{2\beta +5}}\le a_0(\beta )\left (1+\frac {a_+(\beta )}{a_-(\beta )}\right )\int _{t}^0 \frac {f(x)dx}{|x+t|^{2\beta +5}}\le
$$
$$\le a_0(\beta )\left (1+\frac {a_+(\beta )}{a_-(\beta )}\right )\frac 1{|t|^{2\beta +5}} \int _t^0f(x)dx=a_0(\beta )\left (1+\frac {a_+(\beta )}{a_-(\beta )}\right )\frac {s(t,\varphi )}{|t|^{2\beta +5}}.
$$
The lower estimate is obvious:
$$
K_0(t,\varphi )\ge a_0(\beta )\int _{t }^0\frac {f(x)dx}{|t+x|^{2\beta +5}}\ge \frac {2^{-(2\beta +5)}a_0(\beta )}{|t|^{2\beta +5}}\int _{t}^0 f(x)dx=
$$
$$
=\frac {2^{-(2\beta +5)}a_0(\beta )}{|t|^{2\beta +5}}s(t,\varphi ).
$$
Lemma 2 is proved.
\end{proof}

To complete the proof of Theorem 1, it is enough to substitute the estimates of Lemma 2 into relation (4).
\end{proof}

\noindent {\it Remark. } We can take
$$
a(\beta )=aa_0(\beta )2^{-(2\beta+5)}\frac {(1+2\beta )^2}{4(2+\beta )},\quad A(\beta )=Aa_0(\beta )\frac {2(2+\beta )^2}{(1+2\beta )}\left (1+\frac {a_+(\beta )}{a_-(\beta )}\right )
$$
as constants in Theorem 1. Constants $a_0(\beta ),a_\pm (\beta )$ are defined in Lemma 2 and $a, A$ are absolute constants in Theorem A.

\section{Localization of the norm in the spaces $G^\alpha$}
\label{Prel2}

The next theorem enables one to localize the integrals over $\mathbb  C\setminus\overline D$, namely, to take them over $\Omega \setminus \overline D$, where $\Omega  $ is an arbitrary
neighbourhood of $\overline D$.

We denote by $B(z,\varepsilon )$ the circle of radius $\varepsilon $ with centre at a point $z$. If $z=0$, then this point will not be indicated. Let $D(\varepsilon )=D+B(\varepsilon )$ and
$$
R(D)= \inf \{ R>0: \quad D\subset \overline {B(R)}\} .
$$

\noindent {\bf Theorem 2}  {\it Let $\gamma \in G^{(\alpha )} (D)$ and  $\alpha \in [0;\frac 32)$. Then
$$
\int \limits _{\mathbb C\setminus \overline D}\vert \gamma ''(\zeta )\vert  ^2{\textrm {dist } }^{2\alpha }(\zeta
)\, dm(\zeta )\le
$$
$$
\le (1+B_0(\alpha ))(1+B(\alpha  ,D))
 \int \limits _{D(\varepsilon )\setminus \overline D }\vert \gamma   ''
(\zeta   )\vert   ^2{\textrm {dist } }^{2\alpha }(\zeta )\, dm(\zeta ),
$$
and
$$
\int \limits _{\mathbb C\setminus \overline D}\vert \gamma ''(\zeta )\vert  ^2{\textrm {dist } }^{2\alpha +1}(\zeta
)\, dm(\zeta )\le
$$
$$
\le (1+5R(D)B_1(\alpha ))(1+5R(D)B(\alpha  ,D))
 \int \limits _{D(\varepsilon )\setminus \overline D }\vert \gamma   ''
(\zeta   )\vert   ^2{\textrm {dist } }^{2\alpha }(\zeta )\, dm(\zeta ),
$$
where
$$
B_0(\alpha )=4^{2\alpha }(4^{(2-\alpha )}-1)^{-1},\quad B_1(\alpha )=4^{2\alpha }(2^{(3-2\alpha )}-1)^{-1},
$$
$$
B(\alpha ,D)=256\frac {(20R)^{2\alpha}(|\partial D|+\pi \varepsilon )^2}{\pi ^2\varepsilon ^{2(\alpha +1)}},
$$
for arbitrary $\varepsilon \in (0;R(D))$. If $\alpha \in [\frac 32;\frac 52)$, then the same estimates are satisfied under the additional condition $\lim _{|z|\longrightarrow \infty }|z||\gamma (z)|=0 $ with constants $B_0, B_1$ replaced by
$$
B'_0(\alpha )=4^{2\alpha }(2^{( 3-\alpha )}-1)^{-1},\quad B'_1(\alpha )=4^{2\alpha }(2^{(5-2\alpha )}-1)^{-1}.
$$
}
\begin{proof}
We will split the proof into two stages. In the first stage (Lemma 3), an upper estimate for the integral of $\vert \gamma ''(\zeta )\vert   $ over the exterior of the circle $\overline {B(4R(D))}$ will be found using the integral over the set $ {B(4R(D))}\setminus \overline D$. Then (Lemma 4) the integral over $ {B(4R(D))}\setminus \overline D$ will be bounded above by the integral over $D( \varepsilon ) \setminus \overline D$.

In what follows, the number $R(D)$ will be abbreviated to $R$.

\noindent {\bf Lemma 3} {\it Let $\gamma \in G^{(\alpha )} D)$ and  $\alpha \in [0;\frac 32)$. Then the relations
$$
\int \limits _{|\zeta |\ge 4R}\vert \gamma ''(\zeta )\vert ^2{\textrm {dist } }^{2\alpha }(\zeta  )\,
dm(\zeta )\le B_0(\alpha ) \int \limits _{B(4R)\setminus \overline D}\vert \gamma ''(\zeta )\vert ^2{\textrm {dist } }^{2\alpha}(\zeta
)\,
dm(\zeta ),
$$
$$
\int \limits _{\mathbb C\setminus \overline D}\vert \gamma ''(\zeta )\vert  ^2{\textrm {dist } }^{2\alpha +1}(\zeta
)\, dm(\zeta )\le 5R B_1(\alpha ) \int \limits _{D(\varepsilon )\setminus \overline D }\vert \gamma   ''
(\zeta   )\vert   ^2{\textrm {dist } }^{2\alpha }(\zeta )\, dm(\zeta )
$$
hold. If $\alpha \in [\frac 32;\frac 52)$, then the same estimates are satisfied under the additional condition $\lim _{|z|\longrightarrow \infty }|z||\gamma (z)|=0 $ with constants $B_0, B_1$ replaced by $B_0', B_1'$.
}
\begin{proof}
We represent the function $\gamma (\zeta )$ in the form of a Laurent series:
$$
\gamma (\zeta  )=\sum  \limits  _{k=0}^\infty  \frac  {\gamma  _k}{\zeta
^{k+1}}, \quad \vert \zeta \vert >R.
$$
By the hypotheses of the lemma, both this series and the series
$$
\gamma ''(\zeta  )=\sum  \limits  _{k=0}^\infty  \frac  {(k+2)(k+1)\gamma  _k}{\zeta
^{k+3}}=\sum  \limits  _{k=0}^\infty  \frac  {\gamma '' _k}{\zeta
^{k+3}}, \quad \vert \zeta \vert >R,
$$
are uniformly convergent on $\mathbb C\setminus  {B(2R)}$. Let's take some number $t\in [0;2)$. If $\zeta $ satisfies the condition $|\zeta |\ge 4R $, then $\textrm {dist } (\zeta )\le |\zeta |+R< 2|\zeta | $, and therefore
$$
\int \limits _{|\zeta |\ge 4R}\vert \gamma ''(\zeta )\vert ^2{\textrm {dist } }^{2t }(\zeta  )\,
dm(\zeta )\le 2^{2t }  \int \limits _{|\zeta |\ge 4R}\vert \gamma ''(\zeta )
\vert ^2|\zeta  |^{2t }\, dm(\zeta ).
$$
Passing to polar coordinates in the right-hand integral and taking into account the orthogonality of the system $e^{ik\varphi }$ with respect to the measure $d\varphi   $ on $[0;2\pi ]$, we
obtain
$$
\int \limits _{|\zeta |\ge 4R}\vert \gamma ''(\zeta )\vert ^2{\textrm {dist } }^{2t}(\zeta  )\,
dm(\zeta )\le 2^{2t +1}\pi \int \limits _{4R}^{\infty }\sum \limits _{k=0}
^\infty \frac {\vert \gamma ''_k\vert ^2r^{2t +1}\, dr}{r^{2(k+3)}}=
$$
\begin{equation}\label{5}
=2^{2t +1}\pi \sum \limits _{k=0}
^\infty \frac {\vert \gamma ''_k\vert ^2 }{2(k+2-t )(4R)^{2(k+2-t)}}.
\end{equation}
Again using the Laurent series, we find a lower estimate for the integral over the annulus $B(4R)\setminus B(2R)$ by taking into account the fact that the inequality
$\textrm {dist } (\zeta )\ge \frac 12|\zeta |$ holds for those $\zeta$ with $|\zeta |\ge 2R$:
$$
\int \limits _{2R\le |\zeta |< 4R}\vert \gamma ''(\zeta )\vert ^2{\textrm {dist } }^{2t }(\zeta  )\,
dm(\zeta )\ge 2^{-2t }\int \limits _{2R\le |\zeta |< 4R}\vert \gamma ''
(\zeta )\vert ^2|\zeta  |^{2t }\, dm(\zeta )=
$$
\begin{equation}\label{6}
= 2^{-2t +1} \pi \sum \limits _{k=0}^\infty
\frac {\vert \gamma ''_k\vert ^2}{2(k+2-t )(4R)^{2(k+2-t )}}(2^{2(k+2-t)}-1).
\end{equation}
Hence,
$$
\int \limits _{2R\le |\zeta |< 4R}\vert \gamma ''(\zeta )\vert ^2{\textrm {dist } }^{2t }(\zeta  )\,
dm(\zeta )\ge
$$
\begin{equation}\label{7}
\ge 2^{-2t +1}\left (2^{2(2-t )}-1\right ) \pi \sum \limits _{k=0}^\infty
\frac {\vert \gamma ''_k\vert ^2}{2(k+2-t )(4R)^{2(k+2-t )}}.
\end{equation}
This, together with (5), implies the estimate
$$
\int \limits _{|\zeta |\ge 4R}\vert \gamma ''(\zeta )\vert ^2{\textrm {dist } }^{2t }(\zeta  )\,
dm(\zeta )\le
$$
$$
\le 2^{4t }\left (2^{2(2-t )}-1\right )^{-1} \int \limits _{2R\le |\zeta |< 4R}\vert \gamma ''(\zeta )\vert ^2{\textrm {dist } }^{2t }(\zeta  )\,
dm(\zeta ).
$$
Assuming $t=\alpha <\frac 32$, we obtain the first estimate in Lemma 3 for $\alpha \in [0;\frac 32)$. Assuming $t=\alpha +\frac 12<2$ and taking into account that $\textrm {dist } (\zeta )\le 5R$ for $|\zeta |\le 4R$, we obtain the second estimate in Lemma 3 for $\alpha \in [0;\frac 32)$.

If $\lim _{|z|\longrightarrow \infty }|z||\gamma (z)|=0$, then $\gamma ''_k=0$  and the sums in relations (5), (6) are taken over  $k\ge 1$. In this case  we obtain the inequality
$$
\int \limits _{2R\le |\zeta |< 4R}\vert \gamma ''(\zeta )\vert ^2{\textrm {dist } }^{2t }(\zeta  )\,
dm(\zeta )\ge
$$
$$
\ge 2^{-2t +1}\left (2^{2(3-t )}-1\right ) \pi \sum \limits _{k=1}^\infty
\frac {\vert \gamma ''_k\vert ^2}{2(k+2-t )(4R)^{2(k+2-t )}},\quad t\in [0;3),
$$
instead of  (7). Hence, and from relation (5) with $\gamma "_0=0$, we have
$$
\int \limits _{|\zeta |\ge 4R}\vert \gamma ''(\zeta )\vert ^2{\textrm {dist } }^{2t }(\zeta  )\,
dm(\zeta )\le
$$
$$
\le 2^{4t }\left (2^{2(3-t )}-1\right )^{-1} \int \limits _{2R\le |\zeta |< 4R}\vert \gamma ''(\zeta )\vert ^2{\textrm {dist } }^{2t }(\zeta  )\,
dm(\zeta ).
$$
Assuming $t=\alpha $, we obtain the first estimate in Lemma 3 for $\alpha \in [\frac 32;\frac 52)$. Assuming $t=\alpha +\frac 12$ and taking into account that $\textrm {dist } (\zeta )\le 5R$ for $|\zeta |\le 4R$, we obtain the second estimate in Lemma 3 for $\alpha \in [\frac 32;\frac 52)$.

Lemma 3 is proved.
\end{proof}

\noindent {\bf Lemma 4} {\it If $\gamma \in G^\alpha (D)$ and $\alpha >0$, then
$$
\int \limits _{B(4R)\setminus \overline D (\varepsilon )}\vert
\gamma   ''(\zeta   )\vert   ^2{\textrm {dist } }^{2\alpha }(\zeta    )\,dm(\zeta    )\le
$$
$$
\le
256\frac {(20R)^{2\alpha}(|\partial D|+\pi \varepsilon )^2}{\pi ^2\varepsilon ^{2(\alpha +1)}}\int \limits _
{D(\varepsilon )\setminus \overline D  }\vert \gamma   ''(\zeta   )\vert   ^2{\textrm {dist } }^{2\alpha }(\zeta    )
\,dm(\zeta    ).
$$}
\begin{proof}
Let $\zeta \notin \overline D(\varepsilon )$. Since
$$
\textrm {dist }   (\zeta  ,D)\le  \textrm {dist }  (\zeta ,D(\varepsilon ))+\varepsilon, \quad |\partial D( \varepsilon )| =|\partial D|+2\pi \varepsilon ,
$$
by Cauchy's formula, the upper estimate
$$
\vert \gamma ''(\zeta )\vert \le \frac  1{2\pi}\left | \int  \limits  _{\partial
D(\varepsilon /2)}\frac {\gamma ''(z)\, dz}{z-\zeta }\right | \le
\frac {|\partial  D(\varepsilon /2)|}{2\pi \textrm {dist } (\zeta ,D(\varepsilon /2))}
\max \limits _{z\in \partial D(\varepsilon /2)}\vert \gamma ''(z)\vert \le
$$
$$
\le \frac {|\partial  D|+\pi \varepsilon}{2\pi (\textrm {dist } (\zeta ,D)-\varepsilon /2)}
\max \limits _{z\in \partial D(\varepsilon /2)}\vert \gamma ''(z)\vert
$$
holds. Clearly, we have $\textrm {dist }  (\zeta  ,D)\ge \varepsilon $  for $\zeta \notin \overline D(\varepsilon )$ and $x/\left (x-\frac \varepsilon 2\right )\le 2$ for $x\ge \varepsilon $.
In addition, $\textrm {dist } (\zeta )\le 5R$ for $\zeta \in B(4R)$. Therefore
\begin{equation}\label{8}
\vert \gamma ''(\zeta )\vert ^2{\textrm {dist } }^{2\alpha }(\zeta ) \le \frac {5^{2\alpha }R^{2(\alpha -1)}(|\partial D|+\pi \varepsilon )^2}
{\pi ^2}  \max \limits _{z\in \partial D(\varepsilon /2)}\vert \gamma ''(z)\vert
^2
\end{equation}
for $\zeta \in B(4R)\setminus \overline D(\varepsilon)$. If $z\in \partial D(\varepsilon /2)$, then the circle $B(z,\varepsilon  /4)$ lies in the domain $D(3\varepsilon /4)\setminus \overline D$.
In addition, if $w\in \partial B(z,\frac \varepsilon 4)$, then $\textrm {dist } (w)\ge \frac \varepsilon 4$. Since $|\gamma  ''(z)|^2$ is a subharmonic function, we obtain the upper estimate
$$
\vert  \gamma  ''(z)\vert  ^2\le  \frac  {16}{\pi  \varepsilon  ^2}\int  \limits
_{B(z,\varepsilon /4)}\vert  \gamma  ''(w)\vert  ^2\, dm(w)\le
$$
$$\le \frac  {16}{\pi  \varepsilon  ^2}
(\sup \limits _{B(z,\varepsilon /4)}{\textrm {dist } }^{-2\alpha }(w))\int  \limits
_{B(z,\varepsilon /4)}\vert  \gamma  ''(w)\vert  ^2{\textrm {dist } }^{2\alpha }(w)\, dm(w)\le
 $$
$$
\le \frac  {4^{2(\alpha +1)}}{\pi  \varepsilon  ^{2(\alpha +1)}} \int  \limits
_{D(\varepsilon) \setminus \overline D}\vert  \gamma  ''(w)\vert  ^2{\textrm {dist } }^{2\alpha }(w)\, dm(w).
$$
Substituting the above result in (8) and integrating over $B(4R)\setminus \overline  D(\varepsilon ) $, we obtain the assertion of the Lemma 4.

Lemma 4 is proved.
\end{proof}

The estimates in Theorem 2 follow directly from the assertions of Lemmas 3 and 4.
\end{proof}

\section{Proof of the main theorem}
\label{Main}

It suffices to prove the main theorem for domains whose boundaries contain no rectilinear segments and no corners. This follows from the fact that the constants
$c(\beta ,D)$ and $C(\beta ,D)$ depend continuously on $D$. Indeed, let us assume that the theorem is true under these additional conditions on the domain $D$. Let $\varepsilon >0$ be an arbitrary number.Let us inscribe in the set $D(\varepsilon )$ a convex polygon containing $D$ and then replace each side of the polygon by an arc of a sufficiently large circle in
such a way that the resulting domain $D'$ remains convex and is contained in the domain $D(2\varepsilon )$. The boundary of $D'$ contains no rectilinear segments but may have
corners. To eliminate these we pass to the domain $D'+B(0,\varepsilon  )$, and thus obtain a domain whose boundary contains no corners or rectilinear segments, contains $D$,
and is contained in $D(3\varepsilon )$. Since $P_\beta (D)\subset P_\beta (D(\varepsilon ))$, we can apply the main theorem in the spaces $P_\beta (D(\varepsilon ))$ for the functions $F\in P_\beta (D)$
and then we can pass to the limit as $\varepsilon \longrightarrow 0$.

For further proof we need some geometric objects. let $s(x,\varphi  )$ denote the area of the part
$$
D_\varphi = \{ w=x+iy:\ f_1(x)< y < f_2(x),\ -R_\varphi <x<0 \}
$$
of $D$ cut off by the line $\textrm {Re } w=x$. The domain $D_\varphi $ results from the transformation
$w=ze^{i\varphi }-h(\varphi )$ of the domain $D$. Hence, $R_\varphi $ is the distance between the support
lines $L(\varphi )$ and $L(\varphi +\pi )$. Let $l(x,\varphi )$ denote the length of the part of the boundary
of $D_\varphi $ cut off by the same line, and $u(x,\varphi  )$  denote the length of the chord cut off by $D_\varphi $
on that line. We write
$$
\sigma (D)=\inf \limits _{\varphi \in [0;2\pi ]}R_ \varphi .
$$
It is clear that $\sigma (D) $ is the "smallest"\ width of the domain.

By Theorem 1, the norm $\Vert F\Vert _{\beta}$ defined in the main theorem for an entire function
$F$ is equivalent to the triple integral
$$
\int \limits _0^{2\pi } \int \limits _{-\infty }^{\infty  }\int  \limits
_{-\infty }^{0}\frac  { |\gamma ''   (e^{-i\varphi    }(h(\varphi
)-(x+iy)))|^2|x|^{2\beta +3}}{s(x,\varphi )}dxdyd\Delta (\varphi ).
$$
We make the change of variables
$$
\zeta =e^{-i\varphi }(h(\varphi )-x-iy),\quad \theta =\varphi
$$
in this integral. Let us discuss the geometric meaning of the old and new variables.

In what follows, $l(\varphi )$ will denote the directed line $\{ te^{i\varphi },\ -\infty <t< \infty \} $. We
choose the anticlockwise direction of motion along the boundary of the domain $D$ and thus define directions along all tangents to the boundary. Let $L(\varphi  )$ denote the
tangent parallel to and having the same direction as the line $l(\frac \pi 2-\varphi  )$. If $te^{-i\varphi  }$ is
the point of intersection of $L(\varphi )$ and $l(-\varphi  )$, then it is easily verified that $h(\varphi )=t$. If the variables $\varphi \in [0; 2\pi)$,
$x<0$ and $y$ are given, then $\zeta $ is the point lying in the plane and having the coordinates $(h(\varphi )-x; -y)$ in the coordinate system formed
by the lines $l(-\varphi )$ (the axis of abcissae) and $l(\frac \pi 2-\varphi  )$ (the axis of ordinates). In this case, the condition
 $x<0$ means that the support line $L(\varphi )$ separates the point $\zeta $ from the domain $D$. Let us determine the ranges of the variables $\theta  $  and  $\zeta  $. The point $\zeta $ obviously lies outside $\overline D$. For a fixed $\zeta \in \mathbb  C\setminus \overline D$, the angle $\theta$ must satisfy the condition that the
support line $L(\theta )$ separates the point $\zeta $ from the domain $D$. We draw two tangent lines from the point $\zeta $ to the boundary of $D$. Let their directions coincide
with the respective directions of the lines $l(\varphi  _1)$ and $l(\varphi  _2)$ and let $0\le \varphi _1 \le \varphi _2$. Then the angle $\theta  $ varies from $\varphi _-(\zeta )=\frac \pi 2-\varphi _2$ to $\varphi _+(\zeta ) = \frac \pi 2-\varphi _1$. The Jacobian of the transformation from the variables $\varphi ,x,y$ to the variables $\zeta ,\theta  $ is identically
equal to 1, and we have $x=h(\theta )- \textrm {Re } \zeta e^{i\theta }$. Thus,
$$
\int \limits _0^{2\pi } \int \limits _{-\infty }^{\infty  }\int  \limits
_{-\infty }^{0}\frac  { |\gamma ''   (e^{-i\varphi    }(h(\varphi
)-\xi))|^2|x|^{2\beta +3}}{s(x,\varphi  )}dxdyd\Delta  (\varphi   )=
$$
$$=
\int   \limits
_{\mathbb C\setminus \overline D}|\gamma ''(\zeta  )|^2\left  (\int  \limits  _{\varphi
_-(\zeta )}^{\varphi _+(\zeta )}\frac {(\textrm {Re } \zeta e^{i\theta  }-h(\theta
))^{2\beta +3}}{s(h(\theta )-\textrm {Re } \zeta  e^{i\theta   },\theta  )}d\Delta  (\theta
)\right )dm(\zeta ).
$$
The inner integral on the right-hand side of above equality will be denoted by $p(\zeta )$. We have proven the following theorem.

\noindent {\bf Theorem 1$'$}  {\it Let $F=B(\gamma )$ be an entire function satisfying the condition
$$
\Vert F\Vert ^2=\int \limits _0^{2\pi }\int \limits _0^\infty \frac
{\vert F(re^{i\varphi })\vert ^2}{K(re^{i\varphi })r^{2\beta }}\, drd\Delta  (\varphi
)<\infty
$$
and let
$$
p(\zeta )=\int \limits _{\varphi _-(\zeta )}^{\varphi _+(\zeta )}
\frac {(\textrm {Re } \zeta e^{i\theta }-h(\theta ))^{2\beta +3}}{s(h(\theta )-\textrm {Re } \zeta e^{i\theta
},\theta )}\, d\Delta (\theta).
$$
Then
$$
 a(\beta )\Vert F\Vert ^2\le \int \limits _{\mathbb C\setminus \overline D}\vert \gamma
^{''}(\zeta )\vert ^2p(\zeta )\, dm(\zeta )\le A(\beta )\Vert F\Vert ^2,
$$
where the constants $a(\beta ), A(\beta )$ depend only on $\beta $ (see the remark to Theorem 1).
}

\noindent {\bf Lemma 5}  {\it Let
$$
p_0(\zeta )=\int \limits _{\varphi _-(\zeta )}^{\varphi _+(\zeta )}
\frac {(\textrm {Re } \zeta e^{i\theta }-h(\theta ))^{2\beta +2}}{u(h(\theta )-\textrm {Re } \zeta e^{i\theta
},\theta )}\, d\Delta (\theta).
$$
Then

\noindent 1. for points $\zeta  $ such that $\textrm {dist } (\zeta  )\le  \sigma  (D)/2  $, we have
$$
\frac 23p_0(\zeta )\le p(\zeta )\le 2p_0(\zeta );
$$

\noindent 2. If $\textrm {dist } (\zeta )> \sigma (D)/2 $, then we denote by $I_0$ the part of $(\varphi _-(\zeta  );\  \varphi  _+(\zeta  ))  $ in which the
condition
$$
\textrm {Re } \zeta e^{i\varphi }-h(\varphi )\ge  \frac {\sigma (D)}2
$$
holds and by $I$ the remaining part of the interval $(\varphi _-(\zeta  );\  \varphi  _+(\zeta  ))  $. Then
$$
p(\zeta )\le \frac {4 \textrm {diam }^2 (D)|\partial D|}{\sigma^2 (D)|D|}{\textrm {dist } }^{2\beta +3}(\zeta )+2\int \limits _{I}
\frac {(\textrm {Re } \zeta e^{i\theta }-h(\theta ))^{2\beta +2}}{u(h(\theta )-\textrm {Re } \zeta e^{i\theta
},\theta )}\, d\Delta (\theta).
$$}

\begin{proof}
1. By the definition of $h(\varphi )$,
\begin{equation}\label{9}
\textrm {Re } \zeta e^{i\varphi }-h(\varphi )=\min \limits _{z\in \overline D}
(\textrm {Re } \zeta e^{i\varphi }-\textrm {Re } z e^{i\varphi })\le \min \limits _{z\in \overline
D}\vert \zeta -z\vert =\textrm {dist } (\zeta ).
\end{equation}
Therefore, for points $\zeta $ such that the condition $\textrm {dist } (\zeta )<  \sigma  (D)/2$ holds and for
all $\varphi \in (\varphi _-(\zeta );\ \varphi _+(\zeta ))$, the expression $s(h(\varphi  )-\textrm {Re } \zeta  e^{i\varphi  },\  \varphi  )$ can be estimated using
Proposition 2 in [3], and this immediately yields the estimates in the first part of the lemma.

2. If the interval $I_0$ turns out to be empty, then for $\theta  \in  I_0$ the inequality in
assertion (2) of Proposition 2 in [3] can be used to obtain
$$
s(h(\theta  )-\textrm {Re } \zeta  e^{i\theta   }, \theta )\ge \frac {\sigma^2(D)|D|
}{4\textrm {diam }^2 (D)} .
$$
In view of (9) and the geometric meaning of the function $\Delta (\theta )$, the second part of Lemma 5 follows.

Lemma 5 is proved.
\end{proof}

Let us describe in geometric terms the integral appearing in the definition of the
function $p_0(\zeta  )$. The interval $(\varphi _-(\zeta );\varphi  _+(\zeta ))$ of integration consists of the angles $\varphi  $
such that $\textrm {Re }  \zeta   e^{i\varphi }-h(\varphi )\ge 0$, that is, the directions determined by $\varphi   $ for which
the support line $L(\varphi )$ separates the point $\zeta  $ from the domain $D$. In this case, the
expression $\textrm {Re } \zeta e^{i\varphi}-h(\varphi )$ is the distance between $\zeta $ and $L(\varphi  )$. If this support line
undergoes parallel translation by a distance of $\textrm {Re } \zeta e^{i\varphi}-h(\varphi )$, then the domain $D$
cuts off the chord (whose length we denoted by $u(h(\varphi )-\textrm {Re } \zeta e^{i\varphi}, \varphi)$) from the translated
line. Finally, the geometric meaning of the function $\Delta (\varphi )$ is that, for $\varphi _1 \ge \varphi _2$, the difference $\Delta (\varphi  _1)-\Delta  (\varphi  _2)$ is equal to the length of the arc on the boundary of $D$ between the points of tangency of the support lines $L(\varphi _1 )$ and $L(\varphi _2 )$.

This description is unrelated to any coordinate system. In what follows, we
choose a coordinate system related to a fixed point $\zeta  $, partially describe the domain $D$
as the epigraph of a convex function $f$, and describe the integral in the definition
of the function $p_0(\zeta  )$ using $f$.

Let the point $\zeta $ serve as the origin of the coordinate system. There is a single
point $z_0$ on the boundary of $D$ such that
$$
\textrm {dist }  (\zeta  )=\inf  \limits  _{z\in  \partial   D}|z-\zeta
|=|z_0-\zeta |.
$$
We direct the axis of ordinates from $\zeta $ to  $z_0$. In this coordinate system, the
domain $D$ is part of the epigraph of a convex function $f(x)$ defined on the interval
$(X_1; X_2)$ with
$$
X_1=h(\pi ) ,\quad X_2=h(0).
$$
Angles are naturally considered in the new coordinate system. The angles of inclination
with respect to the axis of abscissae for the two tangents to the domain $D$
that pass through the origin will be denoted by $\varphi _1$ and $\varphi  _2$, where $\varphi _1
\le \varphi _2$. Then the upper and lower limits of integration in the integral appearing in the definition
of $p_0$ are $\varphi _-=\frac \pi 2-\varphi _2$ and $\varphi _+=\frac \pi 2-\varphi _1$ respectively. The distance between
the point $\zeta $ and the domain $D$ in this coordinate system is expressed as $f(0)$ or $-h(\pi /2)$.

Furthermore, let us assume that the boundary of $D$ contains no corners or rectilinear
segments (see the remark in the beginning of this section). For the function $f$,
this means that its derivative $f'$ is a strictly increasing continuous function.

As $\theta $ varies from $\varphi _-$ to $\varphi  _+$, the expression $\frac \pi 2-\theta $ varies monotonically from  $\varphi  _2$
to $\varphi _1$, that is, from the angle of inclination of the tangent to the graph of the
function $f$ at the point $X_1$ to that of the tangent at the point $X_2$. Consequently,
if $x(\theta )$ is given by the relation
\begin{equation}\label{10}
f'(x(\theta ))=\tan \left (\frac \pi 2-\theta \right )=\cot (\theta ),
\end{equation}
then the point corresponding to $x(\theta )$ will move monotonically from $X_1$ to $X_2$, and
$(x(\theta ); \ f(x(\theta  )))$ is the support point of the support line $L(\theta )$. The geometrical
meaning of $\Delta (\theta )$ implies that
$$
d\Delta (\theta )=d\left (\int \limits ^{x(\theta )}\sqrt  {1+f'(s)^2}\,
ds\right )=\sqrt  {1+f'(x(\theta ))^2}\, dx(\theta )=
\frac 1{|\sin \theta  |}\,
dx(\theta ).
$$
We have the following representation for the function $p_0$:
\begin{equation}\label{11}
p_0(\zeta  )=\int  \limits  _{\varphi  _-}^{\varphi  _+}\frac  {|h(\theta
)|^{2\beta +2}}{u(h(\theta ),\theta )}\frac 1{|\sin \theta |}\, dx(\theta ),
\end{equation}
where the expressions involved in the integral on the right-hand side are written in
the coordinate system related to the point $\zeta $.

Next, it is necessary to represent $u(h(\theta  ),\theta  )$ using $f$. For this, we introduce some
new functions. We take an arbitrary point $x_0\in [X_1;\ X_2]$ and a positive number $\delta $
and set
$$
\rho   _+(f,x_0,\delta   )=\sup   \left \{ \rho :\   \rho \le    X_2-x_0 , \    \int
 \limits _0^\rho (f'(x_0+y)-f'(x_0))\, dy\le \delta \right \} ,
$$
$$
\rho   _-(f,x_0,\delta   )=\sup   \left \{ \rho :\   \rho \le x_0-X_1, \
\int    \limits _0^\rho (f'(x_0)-f'(x_0-y))\, dy\le \delta \right \} .
$$
$$
\widetilde \rho (f,x_0 ,\delta )= \rho _-(f,x_0 ,\delta )+\rho _+(f,x_0 ,\delta ).
$$

The expression
$$
g(t)=\sup _{x\in [X_1;X_2]}(xt-f(x)),\quad -\infty <t<\infty ,
$$
is called the Young conjugate of $f$. Let $T_1=f'(X_1)$ and $T_2=f'(X_2)$. If $t\in (T_1;T_2)$,
then the supremum in the definition of $g$ is attained at the single stationary point
$x=x(t)$ determined by the condition $f'(x)=t$, that is, $g(t)\equiv x(t)t-f(x(t)),\ t\in (T_1;T_2)$. We differentiate this identity with respect to $t$ and obtain $g'(t)\equiv x(t)$ or
 $g'(f'(x))\equiv x$. Assuming $x=x(\theta )$, from (10) we have
\begin{equation}\label{12}
x(\theta )=g'(\cot (\theta )).
\end{equation}

For a convex function $g$ on $\mathbb R$ and an arbitrary positive number $\delta $, we define the
expression $\rho =\rho (g,t_0,\delta ) $ by the condition
$$
\rho = \sup \left \{ s>0: \ \int \limits  _{-s}^s\vert  g'(t_0+t)-g'(t_0)\vert
\, dt\le \delta \right \} .
$$

From Lemma 5 and from (11), (12) in this paper and from Lemmas 3, 4 in  [3] we obtain the following statement.

\noindent {\bf Lemma 6}  {\it For points $\zeta  $ such that  $\textrm {dist } (\zeta  )\le  \sigma  (D)/2  $, the inequalities
$$
\frac 29\int \limits _{\varphi _-}^{\varphi _+}|h(\theta  )|^{2\beta +1}|\sin  \theta
|\rho \left (g, \cot \theta , \frac {|h(\theta )|}{|\sin \theta |}\right )\, dg' (\cot
\theta) \le p(\zeta )\le
$$
$$\le   \frac {48 \textrm {diam } ^4(D)}{|D|^2}\int
\limits _{\varphi _-}^{\varphi _+}|h(\theta  )|^{2\beta +1}|\sin  \theta
|\rho \left (g, \cot \theta , \frac {|h(\theta )|}{|\sin \theta |}\right )\, dg' (\cot
\theta)
$$
hold.

If $\textrm {dist } (\zeta )> \sigma (D)/2 $, then we denote by $I_0$ the part of $(\varphi _-;\  \varphi  _+)  $ where the condition
$$
-h(\theta )> \frac {\sigma (D)}2
$$
holds and let $I$ be the remaining part of the interval $(\varphi _-;\  \varphi  _+)  $. Then
$$
p(\zeta )\le \frac {4\textrm {diam }^2 (D)|\partial D|}{\sigma^2 (D)|D| }{\textrm {dist } }^{2\beta +3}(\zeta )
+
$$
$$+\frac {48 \textrm {diam } ^4(D)}{|D|^2}\int
\limits _I|h(\theta  )|^{2\beta +1}|\sin  \theta |\rho \left (g, \cot \theta , \frac {|h(\theta )|}{|\sin \theta
|}\right )\, dg' (\cot
\theta).
$$
}
We make the change of variable $\theta =\frac \pi 2-\varphi$ in the integrals and set
\begin{equation}\label{13}
p_1(\zeta )=\int \limits _{\theta _-}^{\theta _+}\left \vert h \left (\frac \pi 2-\varphi
\right )\right \vert ^{2\beta +1}\vert\cos \varphi \vert \rho \left (g, \tan \varphi ,\frac {|h(\frac  \pi  2-\varphi
)|}{\cos \varphi}\right )\, dg'(\tan \varphi),
\end{equation}
where $\theta _\pm =\frac \pi 2-\varphi$ are the angles of inclination with respect to the axis of abscissae of the tangents to the
graph of $f(x)$ that pass through the origin. In the above coordinate system, the distance $\textrm {dist } (\zeta )$ is equal to $-h(\frac  \pi  2)  $, and we
will denote it by $d$. Let us define a point $x=x(\varphi ) \in [X_1;\ X_2]$ using the condition
$$
f'(x(\varphi ))=\tan \varphi .
$$
Then the supremum $\sup _x(xt-f(x)) $ is attained at $x(\varphi )$ for $t=\tan \varphi $, that is,
\begin{equation}\label{14}
g(\tan \varphi )=\tan \varphi \cdot x(\varphi )-f(x(\varphi )).
\end{equation}
Besides, the support line $L(\varphi )$ of $D$ is the tangent to the graph of $f(x)$ at the
point $x(\varphi )$, and $-h(\frac \pi 2-\varphi  )$ is the distance between this tangent line and the origin. It can easily be seen that, in this case, $-h(\frac \pi 2-\varphi  )/ \cos \varphi $ is equal to the ordinate of
the point of intersection of the support line with the axis of ordinates. The equation
of the tangent line through the point $x(\varphi )$ has the form
$$
y=\tan \varphi \cdot (x-x(\varphi ))+f(x(\varphi )),
$$
and therefore
$$
\frac {-h(\frac \pi 2-\varphi  )}{  \cos  \varphi  }=f(x(\varphi  ))-\tan
\varphi\cdot x(\varphi ).
$$
From (14), we obtain that
$$
\frac {-h(\frac \pi 2-\varphi  )}{  \cos  \varphi  }=-g(\tan \varphi ).
$$
Thus, if $\varphi $ monotonically increases from $0$ to $\theta _+$ or decreases from $0$ to $\theta _-$,
then $-g(\tan \varphi )$ monotonically decreases from $d$ to $0$. We set $\varphi _0=0 $ and find the
angles $\varphi _n$ from the conditions
$$
\frac {-h(\frac \pi 2-\varphi _n )}{  \cos  \varphi _n  }=\frac {d}{2^{|n|}}
$$
or, which is the same, from the relations $-g(\tan \varphi _n )=2^{-|n|}d$.
The interval of integration splits into the subintervals $(\varphi _n;\ \varphi _{n+1}] $, $n\in \mathbb Z$, and
the integral in (13) itself can thus be represented as the sum of integrals over these
subintervals. We note that $\rho (g, t, \delta ) $ in a non-decreasing function of the argument  $\delta  $,
and therefore
$$
\rho (g, \tan \varphi ,  2^{-n-1}d)\le \rho \left (g, \tan \varphi ,
\frac {-h(\frac \pi 2-\varphi  )}{  \cos  \varphi }\right )\le
\rho (g, \tan \varphi ,  2^{-n}d),\quad n\ge 0,
$$
$$
\rho (g, \tan \varphi ,  2^{-|n|}d)\le \rho \left (g, \tan \varphi ,
\frac {-h(\frac \pi 2-\varphi  )}{  \cos  \varphi }\right )\le
\rho (g, \tan \varphi ,  2^{-|n|+1}d),\quad n< 0 ,
$$
for $\varphi \in (\varphi _n;\ \varphi _{n+1}]$.  Hence, the change of variables $t=\tan \varphi $ results in
$$
p_1(\zeta )\le \sum \limits _{n=0}^{\infty } \left (\frac d{2^{n}}\right )^{2\beta +1}\int  \limits
_{\tan \varphi _n}^{\tan \varphi _{n+1}}\rho \left (g, t,\frac d{2^n}\right )\, dg'(t)+
$$
\begin{equation}\label{15}
+\sum \limits _{n=-1}^{-\infty }\left (\frac  d{2^{|n|-1}}\right )^{2\beta +1}\int  \limits
_{\tan \varphi _{n}}^{\tan \varphi _{n+1}}\rho \left (g, t,\frac d{2^{|n|-1}}\right )\,
dg'(t),
\end{equation}
\begin{equation}\label{16}
p_1(\zeta ) \ge \left (\frac d2\right )^{2\beta +1}\int \limits _{\tan \varphi _{-1}}^{\tan \varphi  _{1}}
 \frac 1{(1+t^2)^{\beta +1}}\rho
\left (g,t,\frac d2\right )\, dg'(t).
\end{equation}

Let us find upper bounds for $p_1(\zeta )$. We set
$$
t_n=\frac {\tan \varphi _n +\tan \varphi _{n+1}}2.
$$
Then, obviously, $$g(t_n)\ge  g(\tan  \varphi  _n)=-\frac d{2^n}$$ for $n\ge  0$ and $$g(t_n)\ge g(\tan  \varphi  _{n+1})=-\frac d{2^{|n|-1}}$$ for $n<0$. Therefore
$$
g(\tan  \varphi  _n)+g(\tan  \varphi  _{n+1})-2 g(t_n)\le  -\frac d{2^{n+1}}
-\frac d{2^{n}}+\frac {2d}{2^n}=\frac d{2^{n+1}}<\frac d{2^n},\quad n\ge
0,$$
$$
g(\tan   \varphi    _n)+g(\tan    \varphi    _{n+1})-2   g(t_n)\le   \frac
d{2^{|n|}}<\frac d{2^{|n|-1}},\quad n< 0.
$$
Comparing these estimates with definition of $\rho  (g,t,\delta  )$, we obtain
$$
\tan \varphi _{n+1}-t_n=t_n-\tan \varphi _{n}<\rho \left (g,t_n,\frac d{2^n}\right ),\quad n\ge 0,
$$
$$
\tan \varphi _{n+1}-t_n=t_n-\tan \varphi _{n}<\rho \left (g,t_n,\frac d{2^{|n|-1}}\right ),\quad n< 0.
$$
We write
$$
\rho _n =  \rho \left (g,t_n,\frac d{2^n}\right ), \quad n\ge 0 ,\quad
\rho _n = \rho \left (g,t_n,\frac d{2^{|n|-1}}\right ), \quad n< 0.
$$
Formula (15) now implies the inequality
$$
p_1(\zeta )\le  \sum  \limits  _0^  {\infty  }\left (\frac  d{2^n}\right )^{2\beta +1}\int  \limits
_{t_n-\rho _n}^{t_n+\rho _n}\rho \left (g,t,\frac d{2^n}\right )\, dg'(t)+
$$
$$
+\sum \limits _{-1}^{-\infty }\left (\frac d{2^{|n|-1}}\right )^{2\beta +1}\int  \limits
_{t_n-\rho _n}^{t_n+\rho _n}\rho \left (g,t,\frac d{2^{|n|-1}}\right )\, dg'(t).
$$
According to item 3 of Lemma 5 in [3] we obtain
$$
p_1(\zeta )\le \sum _0^\infty \left (\frac d{2^n}\right )^{2\beta +1}\cdot \frac {4d}{2^n}+
\sum _{-1}^{-\infty }\left (\frac d{2^{|n|-1}}\right )^{2\beta +1}\cdot \frac {4d}{2^{|n|-1}}=2\frac {4^{\beta +2}}{4^{\beta +1}-1}
{\textrm {dist } }^{2\beta +2}(\zeta ).
$$
Let $I'=\{ \frac  \pi  2-\theta,\  \theta  \in  I\}  $, where $I\subset  (\varphi  _-;  \varphi  _+)$ is the interval defined in Lemma 6.
Then $I'\subset (\theta _-; \theta _+)$ and
$$
\int \limits _I\vert h(\theta )\vert ^{2\beta +1}\vert \sin \theta \vert \rho \left (g, \cot \theta  ,
\frac {|h(\theta )| }{|\sin \theta | }\right )\, dg'(\cot \theta)=
$$
$$
=\int  \limits
_{I'}\left \vert h\left (\frac \pi  2-\varphi  \right )\right \vert ^{2\beta +1} \vert \cos  \varphi  \vert  \rho  \left (g,  \tan
\varphi  , \frac {|h(\frac \pi 2-\varphi  )| }{|\cos \varphi | }\right )\, dg'(\tan
\varphi)\le p_1(\zeta)\le
$$
$$
\le 2\frac {4^{\beta +2}}{4^{\beta +1}-1}
{\textrm {dist } }^{2\beta +2}(\zeta ).
$$
Substituting the above two inequalities in the relations of Lemma 6, we find upper
estimates for $p(\zeta )$: if $ \textrm {dist } (\zeta )\le \frac {\sigma (D)}2$, then
\begin{equation}\label{17}
 p(\zeta )\le \frac {48\textrm {diam } ^4(D)}{|D|^2}p_1(\zeta)\le \frac {6\cdot 4^{\beta +4}\textrm {diam } ^4(D)}{(4^{\beta +1}-1)|D|^2}{\textrm {dist } }^{2\beta +2}(\zeta
),
\end{equation}
and if $\textrm {dist }  (\zeta  )>  \frac  {\sigma(D)}2$, then
\begin{equation}\label{18}
p(\zeta )\le \frac {4\textrm {diam }^2 (D)|\partial D|}{\sigma^2 (D)|D| }{\textrm {dist } }^{2\beta +3}(\zeta )+\frac {6\cdot 4^{\beta +4}\textrm {diam }
 ^4(D)}{(4^{\beta +1}-1)|D|^2}{\textrm {dist } }^{2\beta +2}(\zeta ).
\end{equation}

Furthermore, using (16), we can find lower bounds for $p_1(\zeta )$ and thus for $p(\zeta  )$. In what follows, we shall assume without loss of generality that
$$
\min (-\tan \  \varphi _{-1}, \tan \ \varphi _1)=\tan \ \varphi _1.
$$
We write $\rho = \rho  \left (g,0,\frac  d2\right )$ for brevity. By the definition of $\rho \left (g,0,\frac  d2\right )$,
$$
g(\rho )+g(-\rho )-2g(0)=\frac d2,
$$
On the other hand,
$$
g(0)=\max \limits _{X_1\le x\le X_2}(-f(x))=-\min \limits _{X_1\le x\le X_2}
f(x)=-f(0)=-d.
$$
Therefore the relation
$$
g(\rho )+g(-\rho )=-\frac 32d
$$
must hold. By the definition of the angles $\varphi _{\pm 1}$, we have
$$
g(\tan \ \varphi _1)+g(-\tan \  \varphi  _1)=-\frac  d2+g(-\tan  \  \varphi
_1)\ge -\frac d2+g(0)=-\frac 32d.
$$
Consequently, the inequality
$$
\tan \ \varphi _1\ge \rho =\rho \left (g,0,\frac d2\right )
$$
holds, which, together with (16), implies that
$$
p_1(\zeta )\ge \left (\frac d2 \right )^{2\beta +1}\int \limits _{-\tan \ \varphi _1}^{\tan \ \varphi _1}
\frac 1{(1+t^2)^{\beta +1}}\rho \left (g, t,\frac d2\right )\, dg'(t)\ge
$$
$$
\ge \left (\frac d 2\right )^{2\beta +1}\frac 1{(1+\tan ^2 \ \varphi _1)^{\beta +1}} \int \limits _{-\rho }^{\rho }\rho \left (g,t,\frac d2\right )\, dg'(t).
$$
Proposition 3 from [3] can now be used to find the upper bound for $\tan  \ \varphi _1$, and the item 3 of Lemma 5 from [3]
can be applied to obtain the lower estimate for the integral:
$$
p_1(\zeta )\ge  \left (\frac  d2 \right )^{2\beta +1} \cdot  \left (1+\frac  {25\textrm {diam }  ^2
(D)}{4\sigma  ^2(D)}\right )^{-(\beta +1)}\cdot \frac d2=
$$
$$
=4^{-(\beta +1)}\left (1+\frac  {25\textrm {diam }  ^2 (D)}{4\sigma  ^2(D)}\right )^{-(\beta +1)}{\textrm {dist } }^{2\beta +2}(\zeta ).
$$
This estimate, together with formulas (17) and (18), enables us to deduce from Lemma 6
the lemma stated below.

\noindent {\bf Lemma 7}  {\it The inequalities
$$
m(\beta ,D){\textrm {dist } }^{2\beta +2}(\zeta )\le p(\zeta ) \le M(\beta ,D)
{\textrm {dist } }^{2\beta +2}(\zeta )
$$
hold for points $\zeta $ such that $\textrm {dist } (\zeta )\le \sigma (D)/2  $.

For points $\zeta $ such that $\textrm {dist } (\zeta )> \sigma (D)/2  $, the estimate
$$
p(\zeta )\le M_0(\beta ,D){\textrm {dist } }^{2\beta +3}(\zeta )+M(\beta ,D){\textrm {dist } }^{2\beta +2}(\zeta )
$$
holds.

Here
$$
m(\beta ,D)=\frac 29\cdot 4^{-(\beta +1)}\left (1+\frac  {25\textrm {diam }  ^2 (D)}{4\sigma  ^2(D)}\right )^{-(\beta +1)},
$$
$$
M(\beta ,D)=\frac {6\cdot 4^{\beta +4} \textrm {diam } ^4\ (D) }{(4^{\beta +1}-1)|D|^2},\quad M_0(\beta ,D)=\frac {4\textrm {diam }^2 (D)|\partial D|}{\sigma^2 (D)|D| },
$$
$\sigma (D)$ is the smallest width of the domain $D$ over all directions, $\textrm {dist } (\zeta ) $
is the distance between the point $\zeta $ and the domain $D$, $|D|$ is the area of $D$, $|\partial D| $ is
the length of the boundary of $D$ and $\textrm {diam } (D)$ is the diameter of $D$.
}

To complete the proof of the main theorem, all we have to do is put together all the proven estimates.

1. Let $\beta \in (-\frac 12;\frac 12)$.

1.1. Let us prove the lower bound in the main theorem. Using Theorem $1'$ and Lemma 7, we have
$$
\|F\|_{P_\beta } ^2 \ge  \frac  1{A(\beta )}\int   \limits   _{\mathbb   C\setminus   \overline   D}|\gamma
''(\zeta)|^2p(\zeta )\, dm(\zeta )\ge
$$
$$
\ge \frac {m(\beta ,D)}{A(\beta )}\int   \limits   _{  D(\sigma (D)/2)\setminus \overline D}|\gamma
''(\zeta)|^2{\textrm {dist } }^{2(\beta +1)}\, dm(\zeta ).
$$
Next we apply the first part of Theorem 2 with $\varepsilon = \sigma (D)/2$ and $\alpha =\beta +1 \in [0;\frac 32)$:
$$
\|F\|_{P_\beta } ^2 \ge   \frac {m(\beta ,D)}{A(\beta )(1+B_0(\beta ,D))(1+B(\beta ,D))}\int   \limits   _{ \mathbb C\setminus \overline D}|\gamma
''(\zeta)|^2{\textrm {dist } }^{2(\beta +1)}\, dm(\zeta ).
$$
Thus, the left inequality in the main theorem in this case is proved with the constant $$c(\beta ,D)=\frac {m(\beta ,D)}{A(\beta )}(1+B_0(\beta ,D ))^{-1}(1+B(\beta ,D))^{-1}.$$

1.2. To prove the upper bound, we apply Theorem $1'$:
$$
\|F\|_{P_\beta } ^2\le \frac 1{a(\beta )}\int \limits _{\mathbb C\setminus \overline D}\vert \gamma
^{''}(\zeta )\vert ^2p(\zeta )\, dm(\zeta ).
$$
Then we apply the second estimate in Lemma 7:
$$
\|F\|_{P_\beta } ^2\le \frac {M(\beta ,D)}{a(\beta )}\int \limits _{\mathbb C\setminus  \overline D}\vert \gamma
^{''}(\zeta )\vert ^2{\textrm {dist } }^{2(\beta +1)}(\zeta )\, dm(\zeta )+
$$
$$+\frac {M_0(\beta ,D)}{a(\beta )}\int \limits _{\mathbb C\setminus  D(\sigma (D)/2)}\vert \gamma
^{''}(\zeta )\vert ^2{\textrm {dist } }^{2\beta +3}(\zeta )\, dm(\zeta ).
$$
Let us estimate the second integral using the second inequality of the first part of Theorem 2 for $\alpha =\beta +1\in [0;\frac 32)$. The right inequality in the main theorem in this case is proved with the constant $$C(\beta ,D)=\frac {M(\beta ,D)}{a(\beta )}+\frac {M_0(\beta ,D)} {a(\beta )}(1+5RB_0(\beta +1 ,D))^{-1}(1+B(\beta +1,D))^{-1}.$$

2. Let $\beta \in [\frac 12; \frac 32)$. If $F\in P_\beta (D)$, then by the definition of $\|F\|_{P_\beta }$, for $\beta \ge \frac 12$ we have
$F(0)=\lim _{|z|\longrightarrow \infty }|z||\gamma (z)| =0$. Therefore we can use the second part of Theorem 2. As a result, estimates of the main theorem will be proved with the same constants as in the case of $\beta <\frac 12$ with the constant $B_0(\beta +1,D)$ replaced by $B'_0(\beta +1,D)$.

\end{document}